\documentclass[journal]{IEEEtran}
\usepackage{graphicx,epsf,psfrag,amsmath,amssymb,amsfonts,verbatim}
\usepackage[mathscr]{eucal}
\usepackage{algorithmic}
\usepackage{verbatim,pifont,color,epsf,psfrag,amssymb,amsfonts,amsmath}
\usepackage[tight,footnotesize]{subfigure}

\newtheorem{theorem}{Theorem}

\def\beq{\begin{equation}}
\def\eeq{\end{equation}}
\def\bea{\begin{eqnarray}}
\def\eea{\end{eqnarray}}
\def\ba{\begin{array}}
\def\ea{\end{array}}
\def\bcenter{\begin{center}}
\def\ecenter{\end{center}}
\def\bitem{\begin{itemize}}
\def\eitem{\end{itemize}}
\def\benum{\begin{enumerate}}
\def\eenum{\end{enumerate}}
\def\bdoc{\begin{document}}
\def\edoc{\end{document}}
\def\defeq{{\stackrel{\Delta}{=}}}

\def\eg{{\it e.g., \/}}
\def\etal{{\it et al. \/}}
\def\viz{{\it viz.,\ \/}}
\def\ie{{\it i.e.,\ \/}}

\definecolor{bgrd}{rgb}{1,1,1}
\definecolor{grey}{rgb}{0.9,0.9,0.6}
\definecolor{gray}{rgb}{0.5,0.5,0.5}
\def\tcr{\textcolor{red}}
\def\tcb{\textcolor{blue}}
\def\tck{\textcolor{black}}
\def\tcy{\textcolor{yellow}}
\def\tcg{\textcolor{green}}
\def\tcbg{\textcolor{bgrd}}
\def\tcm{\textcolor{magenta}}
\def\tcw{\textcolor{white}}
\def\tcgr{\textcolor{gray}}

\renewcommand{\thefootnote}{\fnsymbol{footnote}}

\makeatletter
\newdimen{\captionwidth}
\long\def\@makecaption#1#2{%
\captionwidth .9\hsize
\vskip 10pt%
\setbox\@tempboxa\hbox{#1: #2}%
  \ifdim \wd\@tempboxa >\captionwidth%
    \setbox\@tempboxa\hbox{#1:\hspace*{.5em}}%
    \hfil\parbox{\captionwidth}{\raggedright\hangindent \wd\@tempboxa%
    \hangafter=1\unhbox\@tempboxa#2}\hfill%
  \else\centerline{\box\@tempboxa}%
  \fi
}
\makeatother
\def\scalefig#1{\epsfxsize #1\textwidth}

\def\rank{\mbox{rank}}
\def\T{\mbox{\small T}}
\def\H{\mbox{\small H}}
\def\span{\mbox{span}}
\def\iid{\stackrel{\mbox{\small i.i.d.}}{\sim}}
\def\mRe{\mbox{Re}}
\def\mIm{\mbox{Im}}

\newcommand{\mbbC}{\mathbb{C}}
\newcommand{\mbbE}{\mathbb{E}}
\newcommand{\mbbR}{\mathbb{R}}
\newcommand{\mbbV}{\mathbb{V}}

\newcommand{\Amsc}{\mathscr{A}}
\newcommand{\Bmsc}{\mathscr{B}}
\newcommand{\Cmsc}{\mathscr{C}}
\newcommand{\Dmsc}{\mathscr{D}}
\newcommand{\Emsc}{\mathscr{E}}
\newcommand{\Fmsc}{\mathscr{F}}
\newcommand{\Gmsc}{\mathscr{G}}
\newcommand{\Hmsc}{\mathscr{H}}
\newcommand{\Imsc}{\mathscr{I}}
\newcommand{\Jmsc}{\mathscr{J}}
\newcommand{\Kmsc}{\mathscr{K}}
\newcommand{\Lmsc}{\mathscr{L}}
\newcommand{\Mmsc}{\mathscr{M}}
\newcommand{\Nmsc}{\mathscr{N}}
\newcommand{\Omsc}{\mathscr{O}}
\newcommand{\Pmsc}{\mathscr{P}}
\newcommand{\Qmsc}{\mathscr{Q}}
\newcommand{\Rmsc}{\mathscr{R}}
\newcommand{\Smsc}{\mathscr{S}}
\newcommand{\Tmsc}{\mathscr{T}}
\newcommand{\Umsc}{\mathscr{U}}
\newcommand{\Vmsc}{\mathscr{V}}
\newcommand{\Wmsc}{\mathscr{W}}
\newcommand{\Xmsc}{\mathscr{X}}
\newcommand{\Ymsc}{\mathscr{Y}}
\newcommand{\Zmsc}{\mathscr{Z}}

\def\alphabf{\hbox{\boldmath$\alpha$\unboldmath}}
\def\betabf{\hbox{\boldmath$\beta$\unboldmath}}
\def\gammabf{\hbox{\boldmath$\gamma$\unboldmath}}
\def\deltabf{\hbox{\boldmath$\delta$\unboldmath}}
\def\epsilonbf{\hbox{\boldmath$\epsilon$\unboldmath}}
\def\zetabf{\hbox{\boldmath$\zeta$\unboldmath}}
\def\etabf{\hbox{\boldmath$\eta$\unboldmath}}
\def\iotabf{\hbox{\boldmath$\iota$\unboldmath}}
\def\kappabf{\hbox{\boldmath$\kappa$\unboldmath}}
\def\lambdabf{\hbox{\boldmath$\lambda$\unboldmath}}
\def\mubf{\hbox{\boldmath$\mu$\unboldmath}}
\def\nubf{\hbox{\boldmath$\nu$\unboldmath}}
\def\xibf{\hbox{\boldmath$\xi$\unboldmath}}
\def\pibf{\hbox{\boldmath$\pi$\unboldmath}}
\def\rhobf{\hbox{\boldmath$\rho$\unboldmath}}
\def\sigmabf{\hbox{\boldmath$\sigma$\unboldmath}}
\def\taubf{\hbox{\boldmath$\tau$\unboldmath}}
\def\upsilonbf{\hbox{\boldmath$\upsilon$\unboldmath}}
\def\phibf{\hbox{\boldmath$\phi$\unboldmath}}
\def\chibf{\hbox{\boldmath$\chi$\unboldmath}}
\def\psibf{\hbox{\boldmath$\psi$\unboldmath}}
\def\omegabf{\hbox{\boldmath$\omega$\unboldmath}}
\def\Sigmabf{\hbox{$\bf \Sigma$}}
\def\Upsilonbf{\hbox{$\bf \Upsilon$}}
\def\Omegabf{\hbox{$\bf \Omega$}}
\def\Deltabf{\hbox{$\bf \Delta$}}
\def\Gammabf{\hbox{$\bf \Gamma$}}
\def\Thetabf{\hbox{$\bf \Theta$}}
\def\Lambdabf{\mbox{$ \bf \Lambda $}}
\def\Xibf{\hbox{\bf$\Xi$}}
\def\Pibf{{\bf \Pi}}
\def\thetabf{{\mbox{\boldmath$\theta$\unboldmath}}}
\def\Upsilonbf{\hbox{\boldmath$\Upsilon$\unboldmath}}
\newcommand{\Phibf}{\mbox{${\bf \Phi}$}}
\newcommand{\Psibf}{\mbox{${\bf \Psi}$}}

\def\abf{{\bf a}}
\def\bbf{{\bf b}}
\def\cbf{{\bf c}}
\def\dbf{{\bf d}}
\def\ebf{{\bf e}}
\def\fbf{{\bf f}}
\def\gbf{{\bf g}}
\def\hbf{{\bf h}}
\def\ibf{{\bf i}}
\def\jbf{{\bf j}}
\def\kbf{{\bf k}}
\def\lbf{{\bf l}}
\def\mbf{{\bf m}}
\def\nbf{{\bf n}}
\def\obf{{\bf o}}
\def\pbf{{\bf p}}
\def\qbf{{\bf q}}
\def\rbf{{\bf r}}
\def\sbf{{\bf s}}
\def\tbf{{\bf t}}
\def\ubf{{\bf u}}
\def\vbf{{\bf v}}
\def\wbf{{\bf w}}
\def\xbf{{\bf x}}
\def\ybf{{\bf y}}
\def\zbf{{\bf z}}
\def\rbf{{\bf r}}
\def\xbf{{\bf x}}
\def\ybf{{\bf y}}
\def\Abf{{\bf A}}
\def\Bbf{{\bf B}}
\def\Cbf{{\bf C}}
\def\Dbf{{\bf D}}
\def\Ebf{{\bf E}}
\def\Fbf{{\bf F}}
\def\Gbf{{\bf G}}
\def\Hbf{{\bf H}}
\def\Ibf{{\bf I}}
\def\Jbf{{\bf J}}
\def\Kbf{{\bf K}}
\def\Lbf{{\bf L}}
\def\Mbf{{\bf M}}
\def\Nbf{{\bf N}}
\def\Obf{{\bf O}}
\def\Pbf{{\bf P}}
\def\Qbf{{\bf Q}}
\def\Rbf{{\bf R}}
\def\Sbf{{\bf S}}
\def\Tbf{{\bf T}}
\def\Ubf{{\bf U}}
\def\Vbf{{\bf V}}
\def\Wbf{{\bf W}}
\def\Xbf{{\bf X}}
\def\Ybf{{\bf Y}}
\def\Zbf{{\bf Z}}
\def\Ac{{\cal A}}
\def\Bc{{\cal B}}
\def\Cc{{\cal C}}
\def\Dc{{\cal D}}
\def\Ec{{\cal E}}
\def\Fc{{\cal F}}
\def\Gc{{\cal G}}
\def\Hc{{\cal H}}
\def\Ic{{\cal I}}
\def\Jc{{\cal J}}
\def\Kc{{\cal K}}
\def\Lc{{\cal L}}
\def\Mc{{\cal M}}
\def\Nc{{\cal N}}
\def\Oc{{\cal O}}
\def\Pc{{\cal P}}
\def\Qc{{\cal Q}}
\def\Rc{{\cal R}}
\def\Sc{{\cal S}}
\def\Tc{{\cal T}}
\def\Uc{{\cal U}}
\def\Vc{{\cal V}}
\def\Wc{{\cal W}}
\def\Xc{{\cal X}}
\def\Yc{{\cal Y}}
\def\Zc{{\cal Z}}

\def\cross{\!  \times  \!}
\def\alphat{\alpha^{(t)}}
\def\alphab{\alpha^{(b)}}
\def\betat{\beta^{(t)}}
\def\betab{\beta^{(b)}}
\def\ATbf{{\bf A_\sT}}
\def\AHbf{{\bf A_\sH}}
\def\sspan{\mbox{\it span}}

\def\x{{\mathbf x}}
\def\L{{\cal L}}

\def\T{\mbox{\tiny T}}

\def\figwidth{.5\textwidth}

\begin{document}
%
\title{Impact of Data Quality on Real-Time Locational Marginal Price \thanks{L. Jia, J. Kim, R. J. Thomas, and L. Tong are with the School of Electrical and Computer Engineering, Cornell University, Ithaca, NY 14853, USA. Email: {\tt (lj92, jk752, rjt1, ltong)@cornell.edu}.  Part of this work was presented at HICSS 2012 and PES General Meeting 2012.
}}

\author{Liyan Jia,~
				Jinsub Kim,
        Robert J. Thomas,~\IEEEmembership{Life Fellow,~IEEE,}
        and~Lang Tong,~\IEEEmembership{Fellow,~IEEE}}%

\markboth{submitted to the  IEEE Transactions on Power Systems}{Jia  \MakeLowercase{\textit{et al.}}: Effects of Data Quality on Market Functions}

\maketitle{\let\thefootnote\relax\footnotetext{This work is supported in part by a grant under the DoE CERTS program, the NSF under Grant CNS-1135844, and a PSERC grant.}}

\begin{abstract}
The problem of characterizing impacts of data quality on real-time locational marginal price (LMP) is considered. Because the real-time LMP is computed from the estimated network topology and system state, bad data that cause errors in topology processing and state estimation affect real-time LMP. It is shown that the power system state space is partitioned into price regions of convex polytopes.  Under different bad data models, the worst case impacts of bad data on real-time LMP are analyzed. Numerical simulations are used to illustrate worst case performance for IEEE-14 and IEEE-118 networks.
\end{abstract}


%
\IEEEpeerreviewmaketitle

\section{Introduction}
\label{sec:intro}
\IEEEPARstart{T}{he} deregulated electricity market has two interconnected components. The day-ahead market determines the locational marginal price (LMP) based on the dual variables of the optimal power flow (OPF) solution \cite{Wu&etal:96JRE,Litvinov&etal:04TPS}, given generator offers, demand forecast, system topology, and security constraints.  The calculation of LMP in the day-ahead market does not depend on the actual system operation. In the real-time market, on the other hand, an ex-post formulation is often used (\eg by PJM and ISO-New England \cite{Zheng&Litvinov:06TPS}) to calculate the real-time LMP by solving an incremental OPF problem. The LMPs in the day-ahead and the real-time markets are combined in the final clearing and settlement processes.

The real-time LMP is a function of data collected by the supervisory control and data acquisition (SCADA) system.  Therefore,  anomalies in data, if undetected, will affect prices in the real-time market. While the control center employs a bad data detector to ``clean'' the real-time measurements, miss detections and false alarms will occur inevitably.  The increasing reliance on the cyber system also comes with the risk that malicious  data may be injected by an adversary to affect system and real-time market operations.  An intelligent adversary can carefully design a data attack to avoid detection by the bad data detector.

Regardless of the source of data errors, it is of significant value to assess potential impacts of data quality on the real-time market, especially when a smart grid may in the future deploy demand response based on real-time LMP.  To this end, we are interested in characterizing the impact of {\em worst case data errors} on the real-time LMP.  The focus on the worst case also reflects the lack of an accurate model of bad data and our desire to include the possibility of data attacks.

\subsection{Summary of Results and Organization}
We aim to characterize the worst effects of data corruption on real-time LMP. By ``worst'', we mean the maximum perturbation of real-time LMP caused by bad or malicious data, when a fixed set of data is subject to corruption.  {The complete characterization of worst data impact, however, is not computationally tractable. Our goal here is to develop an optimization based approach to search for {\em locally worst data} by restricting the network congestion to a set of lines prone to congestion.  We then apply computationally tractable (greedy search) algorithms to find the worst data and evaluate the effects of worst data by simulations.}

\begin{figure}[thb]
\centering
\psfrag{C0}[c]{\small ${\Xmsc}_0$}
\psfrag{C1}[c]{\small ${\Xmsc}_1$}
\psfrag{C2}[c]{\small ${\Xmsc}_2$}
\psfrag{C3}[c]{\small ${\Xmsc}_3$}
\psfrag{C4}[c]{\small ${\Xmsc}_4$}
\psfrag{p}[c]{\small $\hat{x}$}
\psfrag{q}[c]{\small $\tilde{x}$}
\includegraphics[width = 0.5\textwidth]{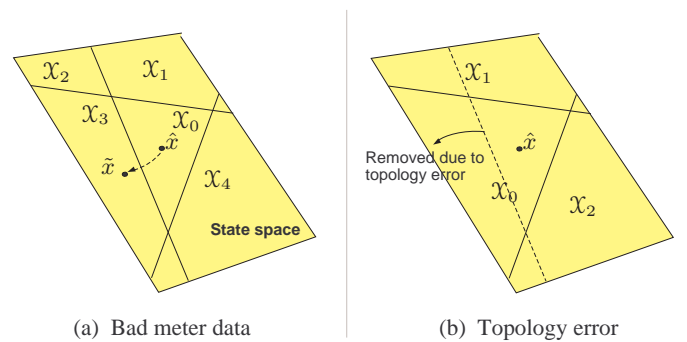}

\caption{Change of real-time LMPs due to bad data.}
\label{fig:geometry_err}
\end{figure}

In characterizing the relation between data and real-time LMP, we first present a geometric characterization of the real-time LMP.  In particular, we show that the state space of the power system is partitioned into polytope price regions, as illustrated in Fig.~\ref{fig:geometry_err}(a), where each polytope is associated with a unique real-time LMP vector, and the price region $\Xmsc_i$ is defined by a particular set of congested lines that determine the boundaries of the price region.

Two types of bad data are considered in this paper.  One is the bad data associated with meter measurements such as the branch power flows in the network. Such bad data will cause errors in state estimation, possibly perturbing, as an example, the correct state estimate $\hat{x}$ in $\Xmsc_0$ to $\tilde{x}$ in $\Xmsc_3$ (as shown in Fig.~\ref{fig:geometry_err}(a)). The analysis of the worst case data then corresponds to finding the worst measurement error such that it perturbs the correct state estimation to the worst price region.



The second type of bad data, one that has not been carefully studied in the context of LMP in the literature, is error in digital measurements such as switch or breaker states.  Such errors lead directly to topology errors therefore causing a change in the polytope structure as illustrated in Fig.~\ref{fig:geometry_err}(b).  In this case, even if the estimated system state changes little, the prices associated with each region change, sometimes quite significantly.

Before characterizing impacts of bad meter data on LMP, we need to construct appropriate models for bad data.  To this end, we propose three increasingly more powerful bad data models based on the dependencies on real-time system measurements: state independent bad data, partially adaptive bad data, and fully adaptive bad data.

In studying the worst case performance, we adopt a widely used approach that casts the problem as one involving an adversary whose goal is to make the system performance as poor as possible.  The approach of finding the worst data is equivalent to finding the optimal strategy of an attacker who tries to perturb the real-time LMP and avoid being detected at the same time.  By giving the adversary more information about the network state and endowing him with the ability to change data, we are able to capture the worst case performance, sometimes exactly and sometimes as bounds on performance.   


Finally, we perform simulation studies using the IEEE-14 and IEEE-118 networks.  We observe that bad data independent of the system state seems to have limited impact on real-time LMPs, and greater price perturbations can be achieved by state dependent bad data.
The results also demonstrate that the real-time LMPs are subject to much larger perturbation if bad topology data are present in addition to bad meter data.

While substantial price changes can be realized for small networks by the worst meter data, as the size of network grows while the measurement redundancy rate remains the same, the influence of worst meter data on LMP is reduced. However, larger system actually gives more possibilities for the  bad topology data to perturb the real-time LMP more significantly.

Our simulation results also show a degree of robustness provided by  the {\em nonlinear state estimator.}   While there have been many studies on data injection attacks based on DC models, very few consider the fact that the control center typically employs the nonlinear WLS state estimator under the AC model.  Our simulation shows that effects of bad analog data designed based on DC model may be mitigated by the nonlinear estimator whereas bad topology data coupled with bad analog data can have greater impacts on LMP.

The rest of the paper is organized as follows. Section~\ref{sec:realtimeLMP} briefly describes a model of real-time LMP and introduces its geometric characterization in the state space of the power system.  
Section~\ref{sec:data_quality} establishes the bad data models and summarizes state estimation and bad data detection procedures at the control center. In Section~\ref{sec:worstcase}, a metric of impact on real-time LMP caused by bad meter data is introduced.  We then discuss the algorithms of finding worst case bad meter data vector in terms of real-time price perturbation under the three different bad data models. Section~\ref{sec:bad_topo} considers the effect of bad topology data on real-time LMP. Finally, in Section~\ref{sec:simulation}, simulation results are presented based on IEEE-14 and IEEE-118 networks.

\subsection{Related Work}
Effects of bad data on power system have been studied extensively in the past,
see \cite{Abur&Exposito:book,%
Handschin&Schweppe&Kohlas&Feichter:75TPAS,Schweppe&Wildes&Rom:70PAS}.  Finding the worst case bad data is naturally connected with the problem of malicious data.  In this context, the results presented in this paper can be viewed as one of analyzing the impact of the worst (malicious) data attack.

In a seminal paper by Liu, Ning, and Reiter \cite{Liu&Ning&Reiter:09CCCS}, the authors first illustrated the possibility that, by compromising enough number of meters, an adversary can perturb the state estimate arbitrarily in some subspace of the state space without being detected by any bad data detector.   Such attacks are referred to as strong attacks.  It was shown by Kosut \etal \cite{Kosut&Jia&Thomas&Tong:11TSG} that the condition for the existence of such undetectable attacks is equivalent to the classical notion of network observability.  

When the adversary can only inject malicious data from a small number of meters,
strong attacks do not exist, and any injected malicious data can be detected with some probability.  Such attacks are referred to as weak attacks \cite{Kosut&Jia&Thomas&Tong:11TSG}.  In order to affect the system operation in some meaningful way, the adversary has to risk being detected by the control center.  The impacts of weak attack on power system are not well understood because the detection of such bad data is probabilistic. Our results are perhaps the first to quantify such impacts. Most related research works focused on DC model and linear estimator while only few have addressed the nonlinearity effect \cite{Jia&Thomas&Tong:12PESGM, HugGiampapa:12TSG}.

It is well recognized that bad data can also cause topology errors \cite{Wu&Liu:89TPS, Clements1988TPS}, and techniques have been developed to detect topology errors.  For instance, the residue vector from state estimation was analyzed for topology error detection \cite{Clements1988TPS, Wu&Liu:89TPS, Costa&Leao:93TPS}.
Monticelli~\cite{Monticelli:93TPS} introduced the idea of generalized state estimation where, roughly speaking, the topology that fits the meter measurements best is chosen as the topology estimate.
The impacts of topology errors on electricity market have not been reported in the literature, and this paper aims to bridge this gap.


The effect of data quality on real-time market was first considered in \cite{Thomas&Tong&Jia&Kosut:10PSerc,Xie&Mao&Sinopoli:10SGC}. In \cite{Xie&Mao&Sinopoli:10SGC}, the authors presented the financial risks induced by the data perturbation and proposed a heuristic technique for finding a case where price change happens. While there are similarities between this paper and \cite{Xie&Mao&Sinopoli:10SGC}, several significant differences exist: (i) This paper focuses on finding the worst case, not only a feasible case. (ii) This paper considers a more general class of  bad data where  bad data may depend dynamically on the actual system measurements rather than static. (iii) This paper considers a broader range of bad data that also include bad topology data, and our evaluations are based on the AC network model and the presence of nonlinear state estimator.

\section{Structures of Real-Time LMP}
\label{sec:realtimeLMP}
In this section, we present first a model for the computation of real-time locational marginal price (LMP).  While ISOs have somewhat different methods of computing real-time LMP, they share the same two-settlement architecture and similar ways of using real-time measurements. In the following, we will use a simplified ex-post real-time market model, adopted by PJM, ISO New England, and other ISOs \cite{Ott:03TPS,Zheng&Litvinov:06TPS}. We view this model as a convenient mathematical abstraction that captures the essential components of the real-time LMP calculation.  For this reason, our results should be interpreted within the specified setup. Our purpose is not to include all details; we aim to capture the essential features.

In real-time, in order to monitor and operate the system, the control center will calculate the estimated system conditions (including bus voltages, branch flows, generation, and demand) based on real-time measurements. We call a branch congested if the estimated flow is larger than or equal to the security limit. The congestion pattern is defined as the set of all congested lines, denoted as $\hat{\mathscr{C}}$.  Note that we use hat (\eg $\hat{\mathscr{C}}$) to denote quantities or sets that are estimated based on real-time measurements.
Details of state estimation and bad data detection are discussed in Section \ref{ssec:se}.

One important usage of state estimation is calculating the real-time LMP. Given 
the estimated congestion pattern $\hat{\mathscr{C}}$, the following linear program is solved to find the incremental OPF dispatch and associated real-time LMP, $\hat{\lambda}=(\hat{\lambda}_i)$ \cite{Ott:03TPS}:
\begin{equation}
\begin{array}{l l}
\mbox{minimize} & \sum c_{i}^\text{\tiny G} \Delta p_i-\sum c_{j}^\text{\tiny L} \Delta d_j\\
\mbox{subjcet to} & \sum \Delta p_{i}=\sum \Delta d_j \\
& \Delta p_i^{\text{min}} \le \Delta p_{i} \le \Delta p_i^{\text{max}}\\
& \Delta d_j^{\text{min}} \le \Delta d_{j} \le \Delta d_j^{\text{max}}\\
& \sum_i A_{ki}\Delta p_i - \sum_j A_{kj} \Delta d_j \le 0,\text{for all }k \in \hat{\mathscr{C}}, \\
\end{array}
\label{eq:real-timeLMP}
\end{equation}
where $\Delta d = (\Delta d_j)$ is the vector of incremental dispatchable load, $\Delta p=(\Delta p_i)$ the vector of incremental generation dispatch, $c^\text{\tiny G} = (c_{i}^\text{\tiny G})$ and $c^\text{\tiny L} = (c_{j}^\text{\tiny L})$ the corresponding real-time marginal cost of generations and dispatchable loads, $\Delta p_i^{\text{min}}$ and $\Delta p_i^{\text{max}}$ the lower and upper bounds for incremental generation dispatch, $\Delta d_i^{\text{min}}$ and $\Delta d_i^{\text{max}}$ the lower and upper bounds for incremental dispatchable load, and $A_{ki}$ the sensitivity of branch flow on branch $k$ with respect to the power injection at bus $i$.

%

%

The real-time LMP at bus $i$ is defined as the overall cost increase when one unit of extra load is added at bus $i$, which is calculated as
\begin{equation}
\hat{\lambda}_i = \eta-\sum_{k \in \hat{\mathscr{C}}}A_{ki}\mu_{k}.
\label{eq:LMPcal}
\end{equation}
where $\eta$ is the dual variable for the load-generation equality constraint, and $\mu_{k}$ is the dual variable corresponding to the line flow constraint in (\ref{eq:real-timeLMP}).

Note that in practice, the control center may use the ex-ante congestion pattern, which is obtained by running a 5 minute ahead security-constrained economic dispatch with the state estimation results and the forecasted loads (for the next five-minute interval) and choosing the lines congested at the dispatch solution \cite{Ott:03TPS, Zhang&Litvinov:06TPS}.  However, to avoid the complication due to ex-ante dispatch calculation, we assume that real-time pricing employs the estimated congestion pattern $\hat{\Cmsc}$ obtained from state estimation results.  By doing so, we attempt to find direct relations among bad data, the state estimate, and real-time LMPs.
Notice that once the congestion pattern $\hat{\mathscr{C}}$ is determined, the whole incremental OPF problem (\ref{eq:real-timeLMP}) no longer depends on the measurement data.
%

Under the DC model, the power system state, $x$, is defined as the vector of voltage phases, except the phase on the reference bus. The power flow vector $f$ is a function of the system state $x$,
\begin{equation}
\centering
f=Fx,
\label{eq:flowequation}
\end{equation}
where $F$ is the sensitivity matrix of branch flows with respect to the system state.

Assume the system has $n+1$ buses. Then, $x \in \mathscr{X}=[-\pi,\pi]^n$, where $\mathscr{X}$ represents the state space. Any system state corresponds to a unique point in $\mathscr{X}$. From (\ref{eq:flowequation}), the branch flow $f$ is determined by the system state $x$. Comparing the flows with the flow limits, we obtain the congestion pattern associated with this state. Hence, each point in the state space corresponds to a particular congestion pattern. 

{We note that the above expression in (\ref{eq:LMPcal}) appears earlier in
\cite{Wu&etal:96JRE} where the role of congestion state in LMP computation was discussed.
 In this paper, our objective is to make explicit the connection between data and LMP.  We therefore need a linkage between data and congestion.  To this end, we note that the power system state, the congestion state, and LMP form a Markov chain, which led to a geometric characterization of LMP on the power system state space, as shown in the following theorem.}

\vspace{3pt}
\begin{theorem}[Price Partition of the State Space]
{Assume that the LMP exists for every possible congestion pattern\footnote{This is equivalent to assuming that the derivative of the optimal value of (\ref{eq:real-timeLMP}) with respect to demand at each bus exists}. Then, the state space $\Xmsc$ is partitioned into a set of polytopes $\{\Xmsc_i\}$ where the interior of each $\Xmsc_i$ is associated with a unique congestion pattern $\Cmsc_i$ and a real-time LMP vector. Each boundary hyperplane of $\Xmsc_i$ is defined by a single transmission line.}
\label{thm:partition}
\end{theorem}

\begin{proof}
For a particular congestion pattern $\Cmsc$ defined by a set of congested lines, the set of states that gives $\Cmsc$ is given by
\[
\Xmsc_i\defeq \{x: F_{i\cdot}x \ge T_i^{\text{max}} \text{ } \forall i \in \Cmsc,
F_{j\cdot}x < T_j^{\text{max}} \text{ } \forall j \notin \Cmsc\},
\]
where $F_{i\cdot}$ is the $i$th row of $F$ (see (\ref{eq:flowequation})), and $T_{j}^{\mbox{max}}$  the flow limit on branch $j$. Since $\Xmsc_i$  is defined by the intersection of a set of half spaces, it is a polytope.

{Given an estimated congestion pattern $\hat{\Cmsc}$, the envelop theorem \cite{MasColell&Whinston:book} implies that for any optimal primal solution and dual solution of (\ref{eq:real-timeLMP}) that satisfy the KKT conditions, (\ref{eq:LMPcal}) always gives the derivative of the optimal objective value with respect to the demand at each bus, which we assume exists, $\ie$ each congestion pattern is associated with a unique real-time LMP vector $\lambda$. Hence, all states with the same congestion pattern share the same real-time LMP, which means each polytope $\Xmsc_i$ in $\Xmsc$ corresponds to a unique real-time LMP vector.}
 
\end{proof}
\vspace{3pt}


Theorem~\ref{thm:partition} characterizes succinctly the relationship between the system state and LMP. As illustrated in Fig.~\ref{fig:geometry_err}(a), if bad data are to alter the LMP in real-time, the size of the bad data has to be sufficiently large so that the state estimate at the control center is moved to a different price region from the true system state.


On the other hand, if some lines are erroneously removed from or added to the correct topology, as illustrated in Fig.~\ref{fig:geometry_err}(b), it affects the LMP calculation in three ways\footnote{In addition to these, the change in topology will affect contingency analysis.  Such effect will appear as changes in contingency constraints in real-time LMP calculation (\ref{eq:real-timeLMP}) \cite{Ott:03TPS}.  However, dealing with contingency constraints will significantly complicate our analysis and possibly obscure the more direct link between bad data and real-time LMP.  Hence, we consider only line congestion constraints in (\ref{eq:real-timeLMP}).}. First, the state estimate is perturbed since the control center employs an incorrect topology in state estimation. Secondly, the price partition of the state space changes due to the errors in topology information.  Third, the shift matrix $A$ in (\ref{eq:real-timeLMP}), which is a function of topology, changes thereby altering prices attached to each price region.



\section{Data  Model and State Estimation}
\label{sec:data_quality}

\subsection{Bad Data  Model}
\label{ssec:data_quality}

\subsubsection{Meter data}
In order to monitor the system, various meter measurements are collected in real time, such as power injections, branch flows, voltage magnitudes, and phasors, denoted by a vector $z\in\mbbR^{m}$. \footnote{Notice here both conventional measurements and PMU measurements can be incorporated. Although PMU data seem to have more direct impact on state estimation and real-time LMP calculation, we won't differentiate the types of measurements in the following discussion.} If there exists bad data $a$ among the measurements, the measurement with bad data, denoted by $z_a$, can be expressed as a function of the system states $x$,
\begin{equation}
z_a = z + a=h(x)+w+a,~~a \in \Amsc,
\label{eq:z=hxwa}
\end{equation}
where $w$ represents the random measurement noise.

We make a distinction here between the measurement noise and bad data; the former accounts for random noise independently distributed across all meters whereas the latter represents the perturbation caused by bad or malicious data. We assume no specific pattern for bad data except that they do not happen everywhere.  We assume that bad data can only happen in a subset of the measurements, $\mathscr{S}$. {We call $\mathscr{S}$ as set of suspectable meters, which means the meter readings with in $\mathscr{S}$ may subject to corruption.} If the cardinality of $\mathscr{S}$ is $k$, the feasible set of bad data $a$ is a $k$-dimensional subspace, denoted as $\Amsc = \{a: a_i = 0\text{ for all } i \notin \mathscr{S}\}$.

We will consider three bad data models with increasing power of affecting state estimates.

M1. {\em State independent bad data}: This type of bad data is independent of real-time measurements.  Such bad data may be the replacement of missing measurements.

M2. {\em Partially adaptive bad data}: This type of bad data may arise from the so-called man in the middle (MiM) attack where an adversary intercepts the meter data and alter the data based on what he has observed. Such bad data can adapt to the system operating state.

M3. {\em Fully adaptive bad data}:  This is the most powerful type of bad data, constructed based on the actual measurement $z=h(x)+w$.  

Note that M3 is in general not realistic.  Our purpose of considering this model is to use it as a conservative proxy to obtain performance bounds for the impact of worst case data.


We assume herein a DC model in which the measurement function $h(\cdot)$ in (\ref{eq:z=hxwa}) is linear. Specifically,
\begin{equation}
\label{eq:DCeq}
z_a=Hx+w+a,~~ a \in \Amsc,
\end{equation}
where $H$ is the measurement matrix. Such a DC model, while widely used in the literature, may only be a crude approximation of the real power system.  By making such a simplifying assumption and acknowledging its weaknesses, we hope to obtain tractable solutions in searching for worst case scenarios.
It is important to note that, although the worst case scenarios are derived from the DC model, we carry out simulations using the actual nonlinear system model.

\subsubsection{Topology data}
Topology data are represented by a binary vector $s\in\{0,1\}^{l}$, where each entry of $s$ represents the state of a line breaker ($0$ for open and $1$ for closed).  The bad topology data is modeled as
\begin{equation}
\label{eq:bad_data_topo}
s_b = s + b~\text{ (mod 2)},~~ b\in\Bmsc,
\end{equation}
where $\Bmsc\subset \{0,1\}^{l}$ is the set of possible bad data.   
  When bad data are present, the topology processor will generate the topology estimate corresponding to $s_b$, and this incorrect topology estimate will be passed to the following operations unless detected by the bad data detector.

\subsection{State Estimation}
\label{ssec:se}

We assume that the control center employs the standard weighted least squares (WLS) state estimator. Under DC model,
\begin{equation}\label{eq:WLS}
\hat{x}=\arg\min_x (z-Hx)^{\tiny{\text{T}}} R^{-1} (z-Hx) = Kz,
\end{equation}
where $R$ is the covariance matrix of measurement noise $w$, and $K \triangleq (H^{\tiny{\text{T}}}R^{-1}H)^{-1}H^{\tiny{\text{T}}}R^{-1}$.

If the noise $w$ is Gaussian, the WLS estimator is also the  maximum likelihood estimate (MLE) of state $x$. By the invariant property of MLE, from (\ref{eq:flowequation}), the maximum likelihood estimate of the branch flows is calculated as
\begin{equation}
\label{eq:WLSflow}
\centering
\hat{f} = F\hat{x}=FKz.
\end{equation}


The congestion pattern used in real-time LMP calculation (\ref{eq:real-timeLMP}) is directly from state estimation and consists of all the estimated branch flows which are larger than or equal to the branch flow limits, $\ie$
\begin{equation}
\label{eq:CongestionPattern}
\centering
\hat{\Cmsc} = \{j: \hat{f}_j \ge T_{j}^{\max} \},
\end{equation}
where $T_{j}^{\max}$ is the flow limit on branch $j$.

In the presence of bad meter data $a$, the meter measurements collected by control center is actually $z_a=Hx+w+a$. By using $z_a$,  
 the WLS state estimate is
\begin{equation}
\hat{x}_a= Kz_a = \hat{x}^* + Ka,
\label{eq:state_move}
\end{equation}
where $\hat{x}^* = Kz$ is the ``correct" state estimate without the presence of the bad data (\ie $a = 0$).

Eq. (\ref{eq:state_move}) shows that the effect of bad data on state estimation is linear.  However, because $a$ is confined in a $k$-dimensional subspace $\Amsc$, the perturbation on the actual system state is limited to a certain direction.

When bad data exist both in meter and topology data, the control center uses a wrong measurement matrix $\bar{H}$, corresponding to the altered topology data, and the altered meter data $z_{a}$.  Then, the WLS state estimate becomes
\begin{equation}
\hat{x}_a= \bar{K}z_a = \bar{K}z + \bar{K}a,
\label{eq:state_topo}
\end{equation}
where $\bar{K} \triangleq (\bar{H}^{\tiny{\text{T}}}R^{-1}\bar{H})^{-1}\bar{H}^{\tiny{\text{T}}}R^{-1}$.  Note that unlike the linear effect of bad meter data, bad topology data affects the state estimate by altering the measurement matrix $H$ to $\bar{H}$.

\subsection{Bad Data Detection}

The control center uses bad data detection to minimize the impact of bad data.  Here, we assume a standard bad data detection used in practice, the $J(\hat{x})$-detector in \cite{Handschin&Schweppe&Kohlas&Feichter:75TPAS}.   In particular, the 
%
$J(\hat{x})$-detector performs the test on the residue error, $r\triangleq z -H\hat{x}$, based on the state estimate $\hat{x}$. From the WLS state estimate (\ref{eq:WLS}), we have
\begin{equation}
\label{eq:residue}
\centering
r = \left(I-H(H^{\tiny{\text{T}}}R^{-1}H)^{-1}H^{\tiny{\text{T}}}R^{-1}\right)z  = Uz.
\end{equation}
where $U \defeq (I-H(H^{\tiny{\text{T}}}R^{-1}H)^{-1}H^{\tiny{\text{T}}}R^{-1})$

The $J(\hat{x})$-detector is a threshold detector defined by
\begin{equation}
r^{\tiny{\text{T}}}R^{-1}r = z^{\tiny{\text{T}}}Wz \begin{array}{c}
\mbox{bad data}\\
\gtrless\\
\mbox{good data}\\
\end{array}
\tau,
\end{equation}
where $\tau$ is the threshold calculated from a prescribed false alarm probability, and $W \defeq U^{\tiny{\text{T}}}R^{-1}U$. 
When the measurement data fail to pass the bad data test, the control center declares the existence of bad data and takes corresponding actions to identify and remove the bad data.  


In this paper, we are interested in those cases when bad data are present while the $J(\hat{x})$-detector fails to detect them.

\section{Impact of Bad Data on LMP}
\label{sec:worstcase}
In this section, we examine the impact of bad data on LMP, assuming that the topology estimate of the network is correct.

{
One thing to notice is that in searching for the ``worst'' case, we take the perspective of the control center, not that of the attacker.  In particular, we look for the worst congestion pattern for the LMP computation, even if this particular congestion pattern is difficult for the attacker to discover.  So the focus here is not how easy it is for an attacker to find a locally worst congestion pattern; it is how much such a congestion pattern affects the LMP.}

\subsection{Average Relative Price Perturbation}
In order to quantify the effect of bad data on real-time price, we need to first define the metric to measure the effect. We define the {\em relative price perturbation} (RPP) as the expected percentage price perturbation caused by bad data. Given that LMP varies at different buses, RPP also varies at different locations.

Let $z_a$ be the data received at the control center and $\lambda_i(z_{a})$ the LMP at bus $i$. The RPP at bus $i$ is a function of bad data $a$, given by
\begin{equation}
\label{eq:metric}
\text{RPP}_i(a)=\mathbb{E}\left(\left|\frac{\lambda_i(z_a)-\lambda_i(z)}{\lambda_i(z)}\right|\right),
\end{equation}
where the expectation is over random state and measurement noise. 

To measure the system-wide price perturbation, we define the {\em average relative price perturbation} (ARPP) by
\begin{equation}
\label{eq:overall metric}
\text{ARPP}(a)=\frac{1}{n+1}\sum_i \text{RPP}_i (a),
\end{equation}
where $n+1$ is the number of buses in the system.

The worst case analysis to be followed can be used for other metrics (e.g., price increase ratios or price decrease ratios, which are closely related to the market participants' gain or loss). Similar results can be showed following the same strategies. However, the comparison among different metrics is beyond the scope of this paper.

\subsection{Worst ARPP under State Independent Bad Data Model}
\label{ssec:constant}

First, we consider the  state independent bad data model (M1) given in Section~\ref{ssec:data_quality}. In this model, the bad data are independent of real-time measurements.


In constructing the state independent worst data, it is useful to incorporate prior information about the state.  To this end, we assume that system state follows a Gaussian distribution with mean  $x_0$, covariance matrix $\Sigma_x$.  Typically, we choose $x_0$ as the day-ahead dispatch since the nominal system state in real-time varies around its day-ahead projection.

In the presence of bad data $a$, the expected state estimate and branch flow estimate on branch $i$ are given by
\begin{equation}
\mathbb{E}[\hat{x}] = x_0 + K a.
\end{equation}
\begin{equation}
\mathbb{E}[f_{i}] = F_{i\cdot} \mathbb{E}[\hat{x}]= F_{i\cdot}x_0 + F_{i\cdot} Ka,
\end{equation}
where $F_{i\cdot}$ is the corresponding row of branch $i$ in $F$.

Our strategy is to make this expected state estimate into the region with the largest price perturbation among all the possible regions, $\hat{\Cmsc}^*$. From (\ref{eq:CongestionPattern}), this means making all the expected branch flows satisfy the boundary condition of $\hat{\Cmsc}^*$,
\begin{equation}
\begin{tabular}{ll}
$\mathbb{E}[f_{i}] \ge T_i^{\text{max}} $ & $\text{for }i \in \hat{\Cmsc}^*$ \\
$\mathbb{E}[f_{i}] \le T_j^{\text{max}} $ & $\text{for }j \notin \hat{\Cmsc}^*$.
\end{tabular}
\label{eq:Ccon}
\end{equation}

However, due to the uncertainty (from both system state $x$ and measurement noise $w$), the actual estimated state after attack, $\hat{x}$, may be different from $\mathbb{E}[\hat{x}]$. 
Therefore, we want to make $\mathbb{E}[\hat{x}]$ at the ``center'' of the desired price region, $\ie$ maximizing the shortest distance from $\mathbb{E}[\hat{x}]$ to the boundaries of the polytope price regions while still holding the boundary constraints. The shortest distance can be calculated as
\begin{equation}
\beta = \mbox{min}\{\tilde{\beta}: |\mathbb{E}[f_{i}] - T^{\text{max}}| \ge \tilde{\beta} \text{ for all i}\}.
\end{equation}

However, the existence of bad data detector prevents the bad data vector $a$ from being arbitrarily large. According to (\ref{eq:residue}), the weighted squared residue with $a$ is
\begin{equation}
r^{\tiny{\text{T}}}R^{-1}r = z_a^{\tiny{\text{T}}} W z_a = (w+a)^{\tiny{\text{T}}} W (w+a).
\end{equation}
since $WHx = 0$

Heuristically, since $w$ has zero mean, the term $a^{\tiny{\text{T}}}Wa$ can be used to quantify the effect of data perturbation on estimation residue. Then we use $a^{\tiny{\text{T}}} W a \le \epsilon$ to control the detection probability in the following optimization.

Therefore, for a specific congestion pattern $\hat{\Cmsc}$, the adversary will solve the following optimization problem to move the state estimate to the ``center'' of the price region $\hat{\Cmsc}$ and keeping the detection probability low.
\begin{equation}
\label{eq:constant}
\begin{tabular}{l l}
$\max_{a\in \Amsc, \tilde{\beta} \ge 0}$ & $\tilde{\beta}$  \\
\text{subject to} & $\mathbb{E}[f_{i}] - \tilde{\beta} \ge T_i^{\text{max}},i \in \hat{\Cmsc}$ \\
& $\mathbb{E}[f_{i}] + \tilde{\beta} < T_j^{\text{max}},j \notin \hat{\Cmsc}$\\
& $a^{\tiny{\text{T}}} W a \le \epsilon$,\\
\end{tabular}
\end{equation}
which is a convex program that can be solved easily in practice. We call a region $\hat{\Cmsc}$ {\em feasible} if it makes problem~(\ref{eq:constant}) feasible.

Among all the feasible congestion patterns, the worst region $\hat{\mathscr{C}}^*$ is chosen as the one giving the largest ARPP.
\begin{equation}
\hat{\mathscr{C}}^*=\mbox{arg } \max_{\hat{\mathscr{C}} \in \Gamma} |\tilde{\lambda}_i-\lambda_i(\hat{\mathscr{C}})|,
\end{equation}
where $\tilde{\lambda}_i$ is the LMP at bus $i$ if the $x_0$ is the system state, and $\Gamma$ the set of all the feasible congestion patterns. Hence, the worst case constant bad data vector is the solution to optimization problem~(\ref{eq:constant}) by setting the congestion pattern as $\hat{\mathscr{C}}^*$.



%

\subsection{Worst ARPP under Partially Adaptive Bad Data}

For bad data model M2, only part of the measurement values in real-time are known to the adversary, denoted as $z_{\text{o}}$. The adversary has to first make an estimation of the system state from the observation and prior distribution, then make the attack decision based on the estimation result.

Without the presence of bad data vector, $\ie$ $a = 0$, the system equation~(\ref{eq:DCeq}) gives
\begin{equation}
z_\text{o} = H_{\text{o}} x + w_{\text{o}},
\end{equation}
where $H_{\text{o}}$ is the rows of $H$ corresponding to the observed measurements and $w_{\text{o}}$ the corresponding part in the measurement noise $w$.

%

The minimum mean square error (MMSE) estimate of $x$ given $z_{\text{o}}$ is given by the conditional mean
\begin{equation}
\mathbb{E}(x|z_{\text{o}})=x_0+\Sigma_x H_{\text{o}}^{\tiny{\text{T}}} (H_{\text{o}} \Sigma_x H_{\text{o}}^{\tiny{\text{T}}})^{-1}(z_{\text{o}}-H_{\text{o}} x_0).
\end{equation}

%

Then, the flow estimate on branch $i$ after attack is
\begin{equation}
\mathbb{E}[f_{i}|z_{\text{o}}] = F_{i\cdot} \mathbb{E}[\hat{x}|z_{\text{o}}].
\end{equation}


Still, we want to move the estimation of state to the ``center''. On the other hand, the expected measurement value $\mathbb{E}[z_a|z_{\text{o}}] = H\mathbb{E}[\hat{z}|z_{\text{o}}] + a$. Again, we need a pre-designed parameter $\epsilon$ to control the detection probability. Therefore, the solution to the following optimization problem is the best attack given congestion pattern $\Amsc$
\begin{equation}
\begin{tabular}{l l}
$\max_{a \in \Amsc,\tilde{\beta} \ge 0}$ & $\tilde{\beta}$  \\
\text{subject to} & $\mathbb{E}[f_{i}|z_{\text{o}}] - \tilde{\beta} \ge T_i^{\text{max}},i \in \hat{\Cmsc}$ \\
& $\mathbb{E}[f_{i}|z_{\text{o}}] + \tilde{\beta} < T_j^{\text{max}},j \notin \hat{\Cmsc}$\\
& $(H\mathbb{E}[z_a|z_{\text{o}}]^{\tiny{\text{T}}}) W (H\mathbb{E}[z_a|z_{\text{o}}]) \le \epsilon$.\\
\end{tabular}
\label{eq:partial}
\end{equation}


This problem is also a convex optimization problem, which can be easily solved. Among all the $\hat{\Cmsc}$'s which make the above problem feasible, we choose the one with the largest price perturbation, denoted as $\hat{\Cmsc}^*$. The solution to problem (\ref{eq:partial}) with $\hat{\Cmsc}^*$ as the congestion pattern is the worst bad data vector.

\subsection{Worst ARPP under Fully Adaptive Bad Data}

Finally, we consider the bad data model M3, in which the whole set of measurements $z$ is known to the adversary. The worst bad data vector depends on the value of $z$. Different from the previous two models, with bad data vector $a$, the estimated state is deterministic without uncertainty. In particular
\begin{equation}
\hat{x} = K z + K a.
\end{equation}
And the estimated flow on branch $i$ after attack is also deterministic
\begin{equation}
\hat{f}_i = F_{i\cdot} \hat{x} =F_{i\cdot} K z + F_{i\cdot} K a.
\end{equation}

%

Similar to the previous two models, congestion pattern is called feasible if there exists some bad data vector $a$ to make the following conditions satisfied:
\begin{equation}
\begin{tabular}{l}
$\hat{f}_i \ge T_i^{\text{max}},i \in \hat{\Cmsc}$ \\
$\hat{f}_i < T_j^{\text{max}},j \notin \hat{\Cmsc}$\\
$(z+a)^{\tiny{\text{T}}}W(z+a) \le \tau, \quad a \in \Amsc$.\\
\end{tabular}
\label{eq:full}
\end{equation}

Among all the feasible congestion patterns, we choose the one with the largest price perturbation, $\hat{\Cmsc}^*$. Any bad data vector $a$ satisfying condition~(\ref{eq:full}) can serve as the worst fully adaptive bad data.

\subsection{A Greedy Heuristic}
\label{ssec:heuristic}
{The strategies presented above are based on the exhaustive search over all possible congestion patterns. Such approaches are not scalable for large networks with a large number of possible congestion patterns.
We now present a greedy heuristic approach aimed at reducing computation cost.  In particular, we develop a
gradient like algorithm that searches among a set of likely congestion patterns.}

{
First, we restrict ourselves to the set of lines that are close to their respective flow limits and look for bad data that will affect the congestion pattern.  The intuition is that it is unlikely that bad data can drive the system state sufficiently far without being detected by the bad data detector.  In practice, the cardinality of such a set is usually very small compared with the systems size.}

{
Second, we search for the worst data locally by changing one line in the congestion pattern at a time.  Specifically,  suppose that a congestion pattern is the current candidate for the worst data.  Given a set of candidate lines that are prone to congestions, we search locally by flipping one line at a time from the congested state to the un-congested state and vice versa.  If no improvement can be made, the algorithm stops.  Otherwise, the algorithm updates the current ``worst congestion pattern'' and continue. The effectiveness of this greedy heuristic is tested in Section \ref{ssec:sim_heuristic}. }

\section{Bad Topology Data on LMP}\label{sec:bad_topo}
So far, we have considered bad data in the analog measurements.
In this section, we include the bad \emph{topology} data, and describe another bad data model.

We represent the network topology by a directed graph $\Gmsc=(\Vmsc,\Emsc)$ where each $i\in \Vmsc$ denotes a bus and each $(i,j)\in\Emsc$ denotes a \emph{connected} transmission line.
For each physical transmission line (\eg a physical line between $i$ and $j$), we assign an arbitrary direction (\eg $(i,j)$) for the line, and $(i,j)$ is in $\Emsc$ if and only if bus $i$ and bus $j$ are connected. 


Bad data may appear in both analog measurements and digital (\eg breaker status) data, as described in Section~\ref{ssec:data_quality}:
\begin{equation}\label{eq:both_attack}
\begin{array}{ll}
z_a &= z + a = (Hx + w) + a,~~~~a\in\Amsc,\\
s_b &= s + b~\text{ (mod 2)},~~~~b\in\Bmsc.
\end{array}
\end{equation}

As in Section~\ref{sec:worstcase}, we employ the adversary model to describe the worst case.
The adversary alters $s$ to $s_{b}$ by adding $b$ from the set of feasible attack vectors $\Bmsc\subset \{0,1\}^{l}$ such that the topology processor produces the ``target'' topology $\bar{\Gmsc}$ as the topology estimate.  
In addition, the adversary modifies $z$ by adding $a\in\Amsc$ 
 such that $z_a$ looks consistent with $\bar{\Gmsc}$.

In this section, we focus on the worst case when the adversary is able to alter the network topology without changing the state estimate\footnote{In general, the adversary can design the worst data to affect both the state estimate and network topology.  It is, however, much more difficult to make such attack undetectable. }.  We also require that such bad data are generated by an adversary causing undetectable topology change, \ie the bad data escape the system bad data detection.
For the worst case analysis, we will maximize the LMP perturbation among the attacks within this specific class.
Even though this approach is suboptimal, the simulation results in Section~\ref{sec:simulation} demonstrate that the resulting LMP perturbation is much greater than the worst case of the bad meter data.

Suppose the adversary wants to mislead the control center with the target topology $\bar{\Gmsc}=(\Vmsc,\bar{\Emsc})$, a topology obtained by \emph{removing}\footnote{ Line addition by the adversary is also possible \cite{Kim&Tong:13ISGT}.  However, compared to line removal attacks, line addition attacks require the adversary to observe a much larger set of meter measurements to design undetectable attacks.  In addition, the number of necessary modifications in breaker data is also much larger: to make a line appear to be connected, the adversary should make all the breakers on the line appear to be closed.  Please see \cite{KimTong:13JSAC} for the detail.}  a set of transmission lines $\Emsc_{\Delta}$ in $\Gmsc$ (\ie $\bar{\Emsc} = \Emsc\setminus\Emsc_{\Delta}$).  We assume that the system with $\bar{\Gmsc}$ is observable: \ie the corresponding measurement matrix $\bar{H}$ has full column rank\footnote{Without observability, the system may not proceed to state estimation and real-time pricing.  Hence, for the adversary to affect pricing, the system with the target topology has to be observable.}.

Suppose that the adversary changes the breaker status such that the target topology $\bar{\Gmsc}=(\Vmsc,\bar{\Emsc})$ is observed at the control center.  Simultaneously, if the adversary introduces bad data $a=\bar{H}x-Hx$, then 
\begin{equation}
z_a =Hx+ a + w = \bar{H}x + w,
\end{equation}
which means that the meter data received at the control center are completely consistent with the model generated from $\bar{\Gmsc}$.  Thus, any bad data detector will not be effective.

It is of course not obvious how to produce the bad data $a$, especially when the adversary can only modify a limited number of measurements, and it may not have access to the entire state vector $x$.  Fortunately, it turns out that $a$ can be generated by observing only a few entries in $z$ without requiring global information (such as the state vector $x$) \cite{Kim&Tong:13ISGT}.

\begin{figure}[!t]
\centering
\psfrag{g}[c]{\small{$\Gmsc$}}
\psfrag{gb}[c]{\small{$\bar{\Gmsc}$}}
\includegraphics[width=0.4\textwidth]{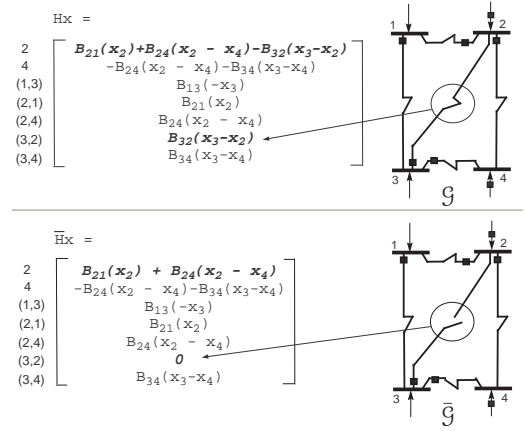}
\caption{$Hx$ and $\bar{H}x$: Each row is marked by the corresponding meter ($i$ for injection at $i$ and $(i,j)$ for flow from $i$ to $j$).}
\label{fig:similarity}
\end{figure}

A key observation is that $Hx$ and $\bar{H}x$ differ only in a few entries corresponding to the modified topology (lines in $\Emsc_\Delta$) as illustrated in Fig.~\ref{fig:similarity}.
 Consider first the noiseless case.
 Let $z_{ij}$ denote the entry of $z$ corresponding to the flow measurement from $i$ to $j$.
As hinted from Fig.~\ref{fig:similarity}, it can be easily seen that $\bar{H}x - Hx$ has the following sparse structure \cite{Kim&Tong:13ISGT}:
\begin{equation}\label{eq:attack_heu}
\bar{H}x - Hx =\displaystyle - \sum_{(i,j)\in \Emsc_{\Delta}}\alpha_{ij}m_{(i, j)},
\end{equation}
where $\alpha_{ij}\in\mbbR$ denotes the line flow from $i$ to $j$ when the line is connected and the system state is $x$, and $m_{(i,j)}$ is the column of the measurement-to-branch incidence matrix, that corresponds to $(i,j)$: \ie $m_{(i,j)}$ is an $m$-dimensional vector with $1$ at the entries corresponding to the flow from $i$ to $j$ and the injection at $i$, and $-1$ at the entries for the flow from $j$ to $i$ and the injection at $j$, and $0$ at all other entries.  Absence of noise implies that $z_{ij} = \alpha_{ij}$, which leads to
\begin{equation}\label{eq:heuristic_attack}
\bar{H}x - Hx =- \sum_{(i,j)\in \Emsc_{\Delta}}z_{ij}m_{(i, j)}.
 \end{equation}

\begin{figure}[t!]
\centering
\psfrag{i}{\small$i$}
\psfrag{j}{\small$j$}
\psfrag{z1}{\small$z_{i}$}
\psfrag{z2}{\small$z_{ij}$}
\psfrag{z3}{\small$z_{ji}$}
\psfrag{z4}{\small$z_{j}$}
\psfrag{z12}{\tcb{\small$z_{i} - z_{ij}$}}
\psfrag{z22}{\tcb{\small$0$}}
\psfrag{z32}{\tcb{\small$0$}}
\psfrag{z42}{\tcb{\small$z_{j}-z_{ji}$}}			
\includegraphics[width = 3.4in]{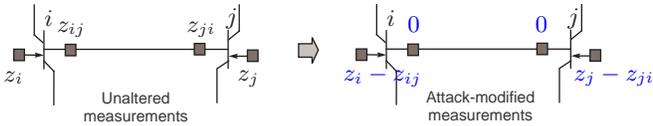}
\caption{
The attack modifies local measurements around the line $(i,j)$ in $\Emsc_{\Delta}$.} 
\label{fig:heuristic}
\end{figure}

With (\ref{eq:heuristic_attack}) in mind, one can see that setting $a=\bar{H}x - Hx$
and adding $a$ to $z$ is equivalent to the following simple procedure:
as described in Fig.~\ref{fig:heuristic}, for each $(i,j)$ in $\Emsc_{\Delta}$,
\begin{enumerate}
\item Subtract $z_{ij}$ and $z_{ji}$ from $z_{i}$ and $z_{j}$ respectively.
\item Set $z_{ij}$ and $z_{ji}$ to be $0$.
\end{enumerate}
where $z_{i}$ is the entry of $z$ corresponding to the injection measurement at bus $i$.

When measurement noise is present (\ie $z = Hx + w$), the idea of the attack is still the same: to make $a$ approximate $\bar{H}x - Hx$ so that $z_{a}$ is close to $\bar{H}x + w$.
Since $z_{ij} = \alpha_{ij} + w_{ij}$, $z_{ij}$ is an unbiased estimate of $\alpha_{ij}$ for each $(i,j)\in\Emsc_{\Delta}$, and this implies that $- \sum_{(i,j)\in \Emsc_{\Delta}}z_{ij}m_{(i, j)}$ is an unbiased estimate of $- \sum_{(i,j)\in \Emsc_{\Delta}}\alpha_{ij}m_{(i, j)} = \bar{H}x - Hx$.  Hence, we set $a$ to be $-\sum_{(i,j)\in \Emsc_{\Delta}}z_{ij}m_{(i, j)}$, the same as in the noiseless setting, and the attack is executed by the same steps as above.

For launching this attack to 
modify the topology estimate from $\Gmsc$ to $\bar{\Gmsc}$, the adversary should be able to (i) set $b$ such that the topology processor produces $\bar{\Gmsc}$ instead of $\Gmsc$ and (ii) observe and modify $z_{ij}$, $z_{ji}$, $z_{i}$, and $z_{j}$ for all $(i,j)\in\Emsc_{\Delta}$.
The attack is feasible if and only if $\Amsc$ and $\Bmsc$ contain the corresponding attack vectors.

To find the worst case LMP perturbation due to undetectable, state-preserving attacks, let $\Fmsc$ denote the set of feasible $\bar{\Gmsc}$s, for which the attack can be launched with $\Amsc$ and $\Bmsc$.
Among the feasible targets in $\Fmsc$, we consider the best target topology that results in the maximum perturbation in real-time LMPs.
If ARPP is used as a metric, the best target is chosen as
\begin{equation}
\bar{\Gmsc}^{*}[z]=\mbox{arg } \max_{\bar{\Gmsc} \in \Fmsc} \sum_{i}\left|\frac{\lambda_i(z;\bar{\Gmsc}) - \lambda_i(z;\Gmsc)}{\lambda_i(z;\Gmsc)}\right|. 
\end{equation}
where $\lambda_{i}(z;\bar{\Gmsc})$ denotes the real-time LMP at bus $i$ when the attack with the target $\bar{\Gmsc}$ is launched on $z$, and $\lambda_i(z;\Gmsc)$ is the real-time LMP under no attack.

\section{Numerical Results}
\label{sec:simulation}

In this section, we demonstrate the impact of bad data on real-time LMPs with the numerical simulations on IEEE-14 and IEEE-118 systems.  We conducted simulations in two different settings: the linear model with the DC state estimator and the nonlinear model with the AC state estimator.  The former is usually employed in the literature for the ease of analysis whereas the latter represents the practical state estimator used in the real-world power system.  
In all simulations, the meter measurements consist of real power injections at all buses and real power flows (both directions) at all branches.  

\subsection{Linear model with DC state estimation}
\label{ssec:sim_linear}

We first present the simulation results for the linear model with the DC state estimator.
We modeled bus voltage magnitudes and phases as Gaussian random variables with the means equal to the day-ahead dispatched values and small standard deviations.
In each Monte Carlo run, we generated a state realization from the statistical model, and the meter measurements were created by the DC model with Gaussian measurement noise.
Once the measurements were created, bad data were added in the manners discussed in Section~\ref{sec:worstcase} and Section~\ref{sec:bad_topo}.  With the corrupted measurements, the control center executed the DC state estimation and the bad data test with the false alarm probability constraint $0.1$.  If the data passed the bad data test, real-time LMPs were evaluated based on the state estimation results.
For IEEE-14 and IEEE-118 system, the network parameters\footnote{
In addition to the network parameters given in \cite{IEEEParameter}, we used the following line limit and real-time offer parameters.  In the IEEE-14 simulation, the generators at the buses 1, 2, 3, 6, and 8 had capacities 330, 140, 100, 100, and 100 MW and the real-time offers 15, 31, 30, 10, and 20 $\$/\text{MW}$.  Lines (2, 3), (4, 5), and (6, 11) had line capacities 50, 50, and 20 MW, and other lines had no line limit.  In the IEEE-118 simulation, the generators had generation costs arbitrarily selected from $\{20, 25, 30, 35, 40 ~\$/\text{MW}\}$ and generation capacities arbitrarily selected from $\{200, 250, 300, 350, 400~\text{MW}\}$.  Total 16 lines had the line capacities arbitrarily selected from $\{70, 90, 110~\text{MW}\}$, and other lines had no line limit.  To handle possible occurrence of price spikes, we set the upper and lower price caps as 500$\$/\text{MW}$ and -100$\$/\text{MW}$ respectively.  Total 1000 Monte Carlo runs were executed for each case.
} are available in \cite{IEEEParameter}.

We used the number of meter data to be modified by the adversary as the metric for the attack effort.  For the 14 bus system, in each Monte Carlo run, we randomly chose two lines, and the adversary was able to modify all the line flow meters on the lines and injection meters located at the ends of the lines.
For the 118 bus system, we randomly chose three lines, and the adversary had control over the associated line and injection meters.  Both state and topology attacks were set to control the same number of meter data\footnote{Topology attacks need to make few additional modifications on breaker state
data such that the target lines appear to be disconnected to the topology processor.  However, for simplicity, we do not take into account this additional effort.} so that we can fairly compare their impacts on real-time LMPs. {As for the meter data attack, we only considered the lines that are close to their flow limits (estimated flows under M1 and M2, or actual flows under M3) as candidates for congestion pattern search. The threshold is chosen as 10MW in our simulation.}

\begin{figure}[!t]
\centering

\subfigure[IEEE-14: ARPP of the worst topology data is $66.1\%$.]
{
\includegraphics[width = 0.4\textwidth]{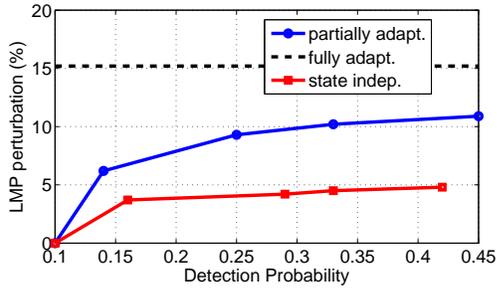}
}

\subfigure[IEEE-118: ARPP of the worst topology data is $22.4\%$.]
{
\includegraphics[width = 0.4\textwidth]{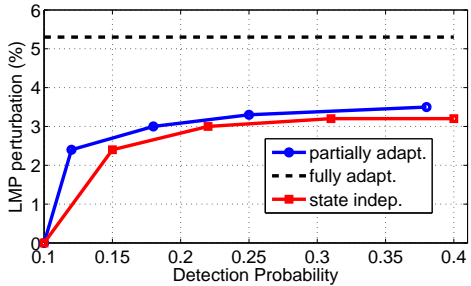}
}
\caption{Linear model: ARPP vs detection prob.}
\label{fig:ARPP_linear}
\end{figure}

Fig.~\ref{fig:ARPP_linear} is the plot of ARPPs\footnote{The detection probabilities for the fully adaptive bad meter data and the bad topology data cases were less than $0.1$ in all the simulations.  In the figures, we draw ARPPs of those cases as horizontal lines so that we can compare them with other cases.
} versus detection probabilities of bad data.  They show that even when bad data were detected with low probability, ARPPs were large, especially for the fully adaptive bad meter data and the bad topology data.

Comparing ARPPs of the three bad meter data models, we observe that the adversary may significantly improve the perturbation amount by exploiting partial or all real-time meter data (for the partially adaptive case, the adversary observed a half of all meters.)  It is worthy to point out that bad topology data result in much greater price perturbation than bad meter data.  


Recall the discussion in Section~\ref{sec:realtimeLMP} and Section~\ref{sec:bad_topo} that bad topology data and bad meter data employ different price-perturbing mechanisms:  bad topology data perturb real-time LMP by restructuring the price regions without perturbing the state estimate (the line-removal attack introduced in Section~\ref{sec:bad_topo} does not perturb state estimate) whereas bad meter data perturb real-time LMP by simply moving the state estimate to a different price region.  Therefore, the observation implies that restructuring the price regions has much greater impact on real-time LMP than merely perturbing the state estimate.

\subsection{Nonlinear model with AC state estimation}

The simulations with the nonlinear model intend to investigate the vulnerability of the real-world power system to the worst adversarial act, designed based on the linear model.
The simulations were conducted on IEEE-14 and IEEE-118 systems in the same manner as the linear case except that we employed the nonlinear model and the AC state estimation.

\begin{figure}[!t]
\begin{center}{
\subfigure[IEEE-14: ARPP of the worst topology data is $95.4\%$.]
{
\includegraphics[width = 0.4\textwidth]{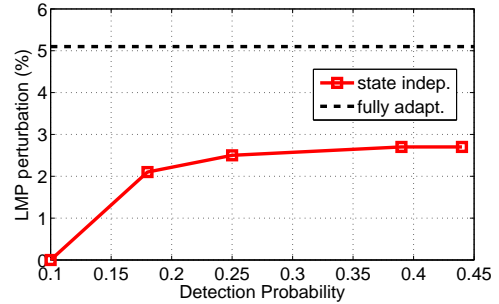}
}

\subfigure[IEEE-118: ARPP of the worst topology data is $76.9\%$.]
{
\includegraphics[width = 0.4\textwidth]{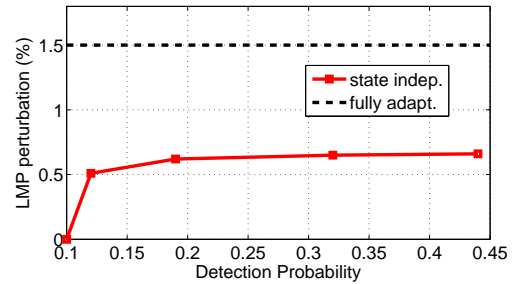}
}
}
\end{center}
\caption{Nonlinear model: ARPP vs detection prob.}
\label{fig:ARPP_nonlinear}
\end{figure}

{Fig.~\ref{fig:ARPP_nonlinear} is the plot of ARPPs versus detection probabilities. The result shows that the proposed methodology can affect the system to some extent even when nonlinear estimator is used, especially when the bad data are present in the topology data, although the nonlinear estimator makes this effect relatively less significant compared with the linear case results. }




\subsection{Performance of the greedy search heuristic}
\label{ssec:sim_heuristic}

{We also conducted simulation based on the proposed greedy search technique in Section \ref{ssec:heuristic}. The simulation was based on 118  bus system, and all parameters were the same as those presented in Section \ref{ssec:sim_linear}. We compared the performance and computation time of the greedy heuristics with exhaustive search benchmark, as shown in Table \ref{tab:greedy}. Notice here the exhaustive search and greedy search are both over the lines that are close to their flow limits (estimated flows under M1 and M2, or actual flows under M3), the same as in Section \ref{ssec:sim_linear}. In Table \ref{tab:greedy}, the second column (average search time) is the average searching time for worst congestion pattern over 1000 Monte Carlo runs, and the third column (accuracy) is the percentage that the greedy search find the same worst congestion pattern as the exhaustive search. From the result, we can see that using greedy heuristic can give us much faster processing algorithm without losing much of the accuracy.}

{
\begin{table}[hbt]
\begin{center}
\caption{Performance of greedy search method}
\begin{tabular}{|c|c|c|}
\hline
method & average search time & accuracy \\
\hline
exhaustive search & 1.23s & - \\
\hline
greedy search & 0.51s & 97.3\% \\
\hline
\end{tabular}
\label{tab:greedy}
\end{center}
\end{table}
}




\section{Conclusion}
We report in this paper a study on impacts of worst data on the real-time market operation.  A key result of this paper is the geometric characterization of real-time LMP given in Theorem \ref{thm:partition}.  This result provides insights into the relation between data and the real-time LMP;  it serves as the basis of characterizing impacts of bad data.

%
%


Our investigation includes bad data scenarios that arise from both analog meter measurements and digital breaker state data.   To this end, we have presented a systematic approach by casting the problem as one involving an adversary injecting malicious data.  While such an approach  often gives overly conservative analysis, it can be used as a measure of assurance when the impacts based on worst case analysis are deemed acceptable.
We note that, because we use adversary attacks as a way to study the worst data, our results have direct implications when cyber-security of smart grid is considered.  Given the increasing reliance on information networks, developing effective countermeasures against malicious data attack on the operations of a future smart grid is crucial.
See \cite{Kosut&Jia&Thomas&Tong:11TSG, Kim&Poor:2011TSG, HugGiampapa:12TSG, KimTong:13JSAC} for discussion about countermeasures.

{From a practical viewpoint, our result can serve as the guideline to the real-time operation. Following the methodology in our paper, worst effect of a specific set of meters on real-time LMP can be checked. Once a huge potential perturbation is detected, alarm should be made and the operator needs to check the accuracy of these specific data, add protection devices, or even add more redundant meters.}

Although our findings are obtained from academic benchmarks involving relatively small size networks, we believe that the general trend that characterizes the effects of bad data is likely to persist in practical networks of much larger size. In particular, as the network size increases and the number of simultaneous appearance of bad data is limited, the effects of the worst meter data on LMP decrease whereas the effects of the worst topology data stay nonnegligible regardless of the network size.  This observation suggests that the bad topology data are potentially more detrimental to the real-time market operation than the bad meter data.

\section*{Acknowledgement}

The authors wish to acknowledge comments and suggestions from the anonymous reviewers that help to clarify a
number of issues and improve the presentation.

\bibliographystyle{IEEEbib}
{
\bibliography{IEEEabrv,ACSP_Bibs/Journal,ACSP_Bibs/Conf,ACSP_Bibs/Misc,ACSP_Bibs/Book}

\begin{thebibliography}{10}
\providecommand{\url}[1]{#1}
\csname url@samestyle\endcsname
\providecommand{\newblock}{\relax}
\providecommand{\bibinfo}[2]{#2}
\providecommand{\BIBentrySTDinterwordspacing}{\spaceskip=0pt\relax}
\providecommand{\BIBentryALTinterwordstretchfactor}{4}
\providecommand{\BIBentryALTinterwordspacing}{\spaceskip=\fontdimen2\font plus
\BIBentryALTinterwordstretchfactor\fontdimen3\font minus
  \fontdimen4\font\relax}
\providecommand{\BIBforeignlanguage}[2]{{%
\expandafter\ifx\csname l@#1\endcsname\relax
\typeout{** WARNING: IEEEtran.bst: No hyphenation pattern has been}%
\typeout{** loaded for the language `#1'. Using the pattern for}%
\typeout{** the default language instead.}%
\else
\language=\csname l@#1\endcsname
\fi
#2}}
\providecommand{\BIBdecl}{\relax}
\BIBdecl

\bibitem{Wu&etal:96JRE}
F.~Wu, P.~Varaiya, P.~Spiller, and O.~S., ``Folk theorems on transimission
  access: proofs and conterexamples,'' \emph{Journal of Regulatory Economics},
  vol.~10, 1996.

\bibitem{Litvinov&etal:04TPS}
E.~Litvinov, T.~Zheng, G.~Rosenwald, and P.~Shamsollahi, ``{Marginal loss
  modeling in LMP calculation},'' \emph{IEEE Transactions on Power Systems},
  vol.~19, no.~2, May 2004.

\bibitem{Zheng&Litvinov:06TPS}
T.~Zheng and E.~Litvinov, ``{Ex-post pricing in the co-optimized energy and
  reserve market},'' \emph{IEEE Transactions on Power Systems}, vol.~21, no.~4,
  November 2006.

\bibitem{Abur&Exposito:book}
A.~Abur and A.~G. Exp\'{o}sito, \emph{Power System State Estimation: Theory and
  Implementation}.\hskip 1em plus 0.5em minus 0.4em\relax CRC, 2000.

\bibitem{Handschin&Schweppe&Kohlas&Feichter:75TPAS}
E.~Handschin, F.~C. Schweppe, J.~Kohlas, and A.~Fiechter, ``{Bad data analysis
  for power system state estimation},'' \emph{IEEE Transactions on Power
  Apparatus and Systems}, vol. PAS-94, no.~2, pp. 329--337, Mar/Apr 1975.

\bibitem{Schweppe&Wildes&Rom:70PAS}
F.~C. Schweppe, J.~Wildes, and D.~P. Rom, ``{Power system static state
  estimation, Parts I, II, III},'' \emph{IEEE Transactions on Power Apparatus
  and Systems}, vol. PAS-89, pp. 120--135, 1970.

\bibitem{Liu&Ning&Reiter:09CCCS}
Y.~Liu, P.~Ning, and M.~K. Reiter, ``False data injection attacks against state
  estimation in electric power grids,'' in \emph{ACM Conference on Computer and
  Communications Security}, 2009, pp. 21--32.

\bibitem{Kosut&Jia&Thomas&Tong:11TSG}
O.~Kosut, L.~Jia, R.~J. Thomas, and L.~Tong, ``Malicious data attacks on the
  smart grid,'' \emph{IEEE Transactions on Smart Grid}, vol.~2, no.~4, pp. 645
  --658, dec. 2011.

\bibitem{Jia&Thomas&Tong:12PESGM}
L.~Jia, R.~J. Thomas, and L.~Tong, ``On the nonlinearity effects on malicious
  data attack on power system,'' in \emph{2012 Power and Energy Society general
  meeting}, July 2012.

\bibitem{HugGiampapa:12TSG}
G.~Hug and J.~Giampapa, ``{Vulnerability assessment of AC state estimation with
  respect to false data injection cyber-attacks},'' \emph{IEEE Transactions on
  Smart Grid}, vol.~3, no.~3, pp. 1362--1370, 2012.

\bibitem{Wu&Liu:89TPS}
F.~F. Wu and W.~E. Liu, ``{Detection of topology errors by state estimation},''
  \emph{IEEE Transactions on Power Systems}, vol.~4, no.~1, pp. 176--183, Feb
  1989.

\bibitem{Clements1988TPS}
K.~Clements and P.~Davis, ``Detection and identification of topology errors in
  electric power systems,'' \emph{IEEE Transactions on Power Systems}, vol.~3,
  no.~4, pp. 1748 --1753, nov 1988.

\bibitem{Costa&Leao:93TPS}
I.~Costa and J.~Leao, ``Identification of topology errors in power system state
  estimation,'' \emph{IEEE Transactions on Power Systems}, vol.~8, no.~4, pp.
  1531 --1538, nov 1993.

\bibitem{Monticelli:93TPS}
A.~Monticelli, ``Modeling circuit breakers in weighted least squares state
  estimation,'' \emph{IEEE Transactions on Power Systems}, vol.~8, no.~3, pp.
  1143 --1149, aug 1993.

\bibitem{Thomas&Tong&Jia&Kosut:10PSerc}
R.~J. Thomas, L.~Tong, L.~Jia, and O.~E. Kosut, ``{Some economic impacts of bad
  and malicious data},'' in \emph{PSerc 2010 Workshop}, vol.~1, Portland Maine,
  July 2010.

\bibitem{Xie&Mao&Sinopoli:10SGC}
L.~Xie, Y.~Mo, and B.~Sinopoli, ``False data injection attacks in electricity
  markets,'' in \emph{Proc. IEEE 2010 SmartGridComm}, Gaithersburg, MD, USA.,
  Oct 2010.

\bibitem{Ott:03TPS}
A.~L. Ott, ``{Experience with PJM market operation, system design, and
  implementation},'' \emph{IEEE Transactions on Power Systems}, vol.~18, no.~2,
  pp. 528--534, May 2003.

\bibitem{Zhang&Litvinov:06TPS}
T.~Zhang and E.~Litvinov, ``{Ex-post pricing in the co-optimized energy and
  reserv markets},'' \emph{IEEE Transactions on Power Systems}, vol.~21, no.~4,
  pp. 1528 -- 1538, Nov. 2006.

\bibitem{MasColell&Whinston:book}
A.~Mas-Colell and M.~D. Whinston, \emph{Microeconomics Theory}.\hskip 1em plus
  0.5em minus 0.4em\relax Oxford University Press, 1995.

\bibitem{Kim&Tong:13ISGT}
J.~Kim and L.~Tong, ``{On topology attack of a smart grid},'' in \emph{2013
  IEEE PES Innovative Smart Grid Technologies (ISGT)}, Washington, DC,
  Feburuary 2013.

\bibitem{KimTong:13JSAC}
------, ``{On topology attack of a smart grid: undetectable attacks and
  countermeasures},'' \emph{IEEE Journal on Selected Areas in Communications},
  vol.~31, no.~7, July 2013.

\bibitem{IEEEParameter}
\BIBentryALTinterwordspacing
``{Power Systems Test Case Archive}.'' [Online]. Available: {\tt
  http://www.ee.washington.edu/research/pstca/}
\BIBentrySTDinterwordspacing

\bibitem{Kim&Poor:2011TSG}
T.~Kim and H.~Poor, ``Strategic protection against data injection attacks on
  power grids,'' \emph{IEEE Transactions on Smart Grid}, vol.~2, no.~2, pp. 326
  --333, june 2011.

\end{thebibliography}
}

\end{document}